\documentclass[12pt]{article}

\usepackage{amssymb}
\usepackage{amsfonts}
\usepackage{amsmath}

\newtheorem{theorem}{Theorem}[section]

\newtheorem{lemma}[theorem]{Lemma}

\newtheorem{proposition}[theorem]{Proposition}

\newtheorem{definition}[theorem]{Definition}
\newtheorem{remark}[theorem]{Remark}

\newcommand{\vol}{V\!ol}

\newcommand{\proof}{\textbf{Proof:\;\;}}
\newcommand{\qed}{\hfill $\blacksquare$\par\medskip}

\begin{document}

\title{Large distortion dimension reduction using random variables }

\author{Alon Dmitriyuk \and Yehoram Gordon
%\thanks{Partially supported by the France-Israel
%Scientific Exchange Program No. 3-4301, and by the Fund for the
%Promotion of Research at the Technion}}
}
\date{}
\maketitle

\begin{abstract}
Consider a random matrix $H:\mathbb{R}^n\longrightarrow\mathbb{R}^m$. Let $D\geq2$ and let $\{W_l\}_{l=1}^{p}$ be a set of $k$-dimensional affine subspaces of $\mathbb{R}^n$. We ask what is the probability that for all $1\leq l\leq p$ and $x,y\in W_l$,
\[
  \|x-y\|_2\leq\|Hx-Hy\|_2\leq D\|x-y\|_2.
\]
We show that for $m=O\big(k+\frac{\ln{p}}{\ln{D}}\big)$ and a variety of different classes of random matrices $H$, which include the class of Gaussian matrices, existence is assured and the probability is very high. The estimate on $m$ is tight in terms of $k,p,D$.
\end{abstract}

\section{Introduction}
The small distance distortion by random linear maps which reduce dimension as applied to sets in Hilbert space is a well known phenomenon. A short history begins with the classical result by Johnson and Lindenstrauss \cite{jl} which stated the dimension reduction scheme for a set of $n$ points in Hilbert space. This was extended by Gordon \cite{gordon88} to a general set and maps and other results concerning finite dimensional Banach spaces using the Gaussian Min-Max theorem proved earlier by the author in \cite{gordon85}. Mendelson, Pajor and Tomczak-Jaegermann proved in \cite{mpj} the dimension reduction scheme using $\psi_2$ random matrices when the range space is Hilbert. As in all these results, including many others, are true for random linear maps with very high probability, these have applications in industry and computer science \cite{achlioptas,av,ct,donoho,im,magen}. We are therefore concerned with the probability estimates as well.

However, the theory on large distance distortion is a different matter and requires different methods in proofs. The paper is concerned with large controlled distance distortion of  preassigned magnitude $D$ by random linear maps $H\in L(\mathbb{R}^n,\mathbb{R}^m)$ where $H$ is a matrix whose rows satisfy certain conditions. To cite one example of large distortion is Theorem 2.8 in \cite{dg2} where $G=(g_{j,i})_{j=1}^{m}\ _{i=1}^{n}$ is a Gaussian matrix with i.i.d. standard Gaussian entries. We recall another of the main results in \cite{dg2} which holds for Gaussian matrices $G$ for large and small distortions:
\begin{theorem}\label{ref}
There is a positive constant $c$ such that the following
holds: Given $0< \varepsilon <\infty$ and $1\leq k\leq n$ and $p$ subspaces
$\{W_l\}_{l=1}^{p}$ of dimension at most $k$ in $\ell_2^n$ and any
$m\geq c\left((1+\varepsilon^{-2})k+\frac{1+\varepsilon}{\varepsilon\ln(1+\varepsilon)}\ln{p}\right)$ and a Gaussian $m\times n$ matrix $G$ with i.i.d. standard Gaussian entries, there is a number $E$ such that the probability that for
every $1\leq l\leq p$ and $x,y\in W_l$
\[
  E\|x-y\|_2\leq\|G(x)-G(y)\|_2\leq E(1+\varepsilon)\|x-y\|_2
\]
is at least $1-\frac{2}{p}$.
\end{theorem}
In \cite{dg2}, for small $\varepsilon$, $0<\varepsilon<1$, we also extended this result to the case when the target is a general Banach space. When the target space was Euclidean, we showed that the corresponding result holds for $\psi_2$ matrices as well, for small distortions.

In this paper we are concerned with large distortions, that is $1\leq \varepsilon$. We show that in this case Theorem \ref{ref} holds for all random matrices with independent $\psi_2$ entries which are absolutely continuous and have jointly bounded density functions. In particular, this class includes Gaussian matrices as well.

We include the case when the rows of the random matrix are independent vectors which are uniformly distributed on $S^{n-1}$.

Note that the results for large distortions contained here and in \cite{dg2} give us lower dimensions at the cost of higher distortions. This is a mathematically intriguing subject, aspects of which were investigated initially by Bourgain \cite{bourgain} and Johnson and Lindenstrauss \cite{jl}. Thus, for example, they showed we can embed an $n$-point metric space into a Hilbert space with distortion $O(\log{n})$ and the dimension of the Hilbert space can be $O(\log{n})$. However, we show here that if we allow for larger distortions, say $O(\log^2{n})$, then using $\varepsilon=\log{n}$ and $k=1$ and $p={n \choose 2}$ in Theorem \ref{ref}, our result shows that we can embed an $n$-point metric space into $O(\frac{\log{n}}{\log\log{n}})$-dimensional Hilbert space, which improves the $O(\log{n})$ above, however, with a larger distortion of $O(\log^2{n})$.

\section{preliminaries}
\begin{definition}
Let $X$ be a random variable. We denote
\[
  F_X=\bigcup_{n=1}^\infty\Big\{\sum_{i=1}^{n}a_i X_i:\sum_{i=1}^{n}a_i^2=1\Big\}
\]
where $X_i\sim X$ and independent. For $\varepsilon>0$ define the $\varepsilon$-concentration of $X$
\[
  C_\varepsilon(X)=\sup_{L\in\mathbb{R}}P(|X-L|<\varepsilon)
\]
and the $\varepsilon$-concentration of $F_X$
\[
  C_\varepsilon(F_X)=\sup_{Y\in F_X}C_\varepsilon(Y).
\]
\end{definition}
\begin{remark}
Note that if $X$ has a density function bounded by $\alpha$ then $C_\varepsilon(X)\leq2\alpha\varepsilon$.
\end{remark}

A result of Ball and Nazarov \cite{bo} states:
\begin{theorem}\label{ob}
There is $c_1>0$ such that for any random variable $X$ and $\varepsilon>0$
\[
  C_\varepsilon(F_X)\leq c_1\cdot C_\varepsilon(X).
\]
\end{theorem}
\medskip

We shall need the following technical lemma later on.
\begin{lemma}\label{lemma01}
Let $X_1,\ldots,X_m$ be independent random variables with density functions which are jointly bounded by $\alpha>0$. Then for all $\lambda>0$
\[
  Pr\left(\sum_{i=1}^{m}X_i^2\leq\lambda m\right)\leq(6\alpha)^m(\lambda)^\frac{m}{2}.
\]
\end{lemma}
\proof For $m=1$ the lemma is clear. For $m\geq2$, using the Stirling approximation for the $\Gamma$ function we obtain for all $\lambda>0$
\begin{eqnarray*}
  &&Pr\left(X_1^2+\cdot+X_m^2\leq\lambda m\right)\\
  &=&\int_{x_1^2+\cdot+x_m^2\leq\lambda m}f_{X_1}(x_1)\cdots f_{X_m}(x_m)dx_1\cdots dx_m\\
  &\leq&\int_{x_1^2+\cdot+x_m^2\leq\lambda m}\alpha^mdx_1\cdots dx_m=
  \alpha^m\vol(\sqrt{\lambda m}B_2^m)=
  \alpha^m\lambda^\frac{m}{2}m^\frac{m}{2}\vol(B_2^m)\\
  &=&\alpha^m\lambda^\frac{m}{2}m^\frac{m}{2}
  \frac{\pi^\frac{m}{2}}{\Gamma\big(\frac{m}{2}+1\big)}\leq
  \alpha^m\lambda^\frac{m}{2}m^\frac{m}{2}\pi^\frac{m}{2}e^{\frac{m}{2}+1} \big(\frac{m}{2}+1\big)^{-\frac{m}{2}}\leq(6\alpha)^m(\lambda)^\frac{m}{2}.
\end{eqnarray*}

\qed

\begin{remark}
Taking $X_i\sim N(0,\sigma^2)$ it is easy to show that Lemma \ref{lemma01} is tight in terms of $\alpha,m,\lambda$.
\end{remark}

\begin{definition}
Let $p\geq1$. The $\psi_p$ norm of a random variable $X$ is
\[
  \|X\|_{\psi_p}=\inf\{C>0:\mathbb{E}\exp\frac{|X|^p}{C^p}\leq2\}.
\]
For a metric space $(T,d)$, an \textit{admissible} sequence of $T$ is a collection of subsets $\{T_s\}_{s\geq0}$ such that for every $s\geq1$ we have $|T_s|\leq2^{2^s}$ and $|T_0|=1$. The $\gamma_p$ functional is defined by
\[
  \gamma_p(T,d)=\inf_{\{T_s\}_{s\geq0}}\sup_{t\in T}\sum_{s=0}^{\infty}2^\frac{s}{p}d(t,T_s)
\]
where the infimum is taken over all admissible sets.
\end{definition}

The upper bound in the following result is by Fernique \cite{fernique} and the lower bound is by Talagrand \cite{Talagrand}.

\begin{theorem}\label{ft}
There are positive constants $c_2,c_3$ such that whenever $\{G_t\}_{t\in T}$ is Gaussian process indexed by a metric space $(T,d)$ for which $d^2(s,t)=\mathbb{E}|G_s-G_t|^2$ then
\[
  c_2\gamma_2(T,d)\leq\mathbb{E}\sup_{t\in T}|G_t|\leq c_3\gamma_2(T,d).
\]
\end{theorem}
We shall use Theorem \ref{ft} in the case where the metric space $T$ is a subset of the unit sphere of $\mathbb{R}^n$ with the Euclidean metric and the Gaussian process is $G_t=\sum_{i=1}^{n}t_ig_i$ where $\{g_i\}_{i=1}^{n}$ are i.i.d. standard Gaussian variables.
\begin{lemma}\label{lemma02}
Let $f:(T,d)\longrightarrow (S,\rho)$ be a one to one correspondence between metric spaces with Lipschitz constant $\beta$. Then $\gamma_p(S,\rho)\leq\beta\gamma_p(T,d)$ for all $p\geq1$.
\end{lemma}
\proof Since $f$ is a one to one correspondence then $\{T_s\}_{s=0}^\infty\longrightarrow\{f(T_s)\}_{s=0}^\infty$ is a one to one correspondence between admissible sequences and therefore
\begin{eqnarray*}
  \gamma_p(S,\rho)=\inf_{\{T_s\}}
  \sup_{t\in T}\sum_{s=0}^{\infty}2^\frac{s}{p}\rho(f(t),f(T_s))\leq
  \inf_{\{T_s\}}
  \sup_{t\in T}\sum_{s=0}^{\infty}2^\frac{s}{p}\beta d(t,T_s)=\beta\gamma_p(T,d).
\end{eqnarray*}
\qed

For completeness we state the next well known result and its proof is contained in the appendix.
\begin{proposition}\label{psi2}
Let $p\geq1$ and let $X$ be a random variable.\\
1. If $\|X\|_{\psi_p}=a$ then $Pr(|X|>t)\leq2\exp\left(-\frac{t^p}{a^p}\right)$ for all $t>0$.\\
2. If $Pr(|X|>t)\leq b\exp\left(-\frac{t^p}{a^p}\right)$ for all $t\geq0$ then $\|X\|_{\psi_p}\leq (b+1)^\frac{1}{p}\cdot a$.\\
3. If $\mathbb{E}\exp(tX)\leq\exp(c^2t^2)$ then $Pr(|X|>t)\leq2\exp\left(-\frac{t^2}{4c^2}\right)$ and therefore $\|X\|_{\psi_2}\leq\sqrt{12}c$.\\
4. If $\|X\|_{\psi_2}=a$ and $\mathbb{E}X=0$ then $\mathbb{E}\exp(tX)\leq\exp(2a^2t^2)$. If $X$ is symmetric then $\mathbb{E}\exp(tX)\leq\exp(a^2t^2)$.\\
5. If $\{X_i\}_{i=1}^{n}$ are independent and satisfy $\mathbb{E}X_i=0$ for $1\leq i\leq n$, and $\max_{1\leq i\leq n}\|X_i\|_{\psi_2}=a<\infty$ then $\|\sum_{i=1}^{n}a_iX_i\|_{\psi_2}\leq4a$ whenever $\sum_{i=1}^{n}a_i^2=1$. If $X_1,\ldots,X_n$ are symmetric then $\|\sum_{i=1}^{n}a_iX_i\|_{\psi_2}\leq\sqrt{12}a$  whenever $\sum_{i=1}^{n}a_i^2=1$.\\
6. If $\{X_i\}_{i=1}^{n}$ are independent and satisfy $\mathbb{E}\exp(tX_i)\leq\exp(c^2t^2)$ for $1\leq i\leq n$, then $\mathbb{E}\exp(t\sum_{i=1}^{n}a_iX_i)\leq\exp(c^2t^2)$ whenever $\sum_{i=1}^{n}a_i^2=1$.\\
7. If $X$ is a random variable for which $|X|\leq a$ then $\|X\|_{\psi_p}\leq4^\frac{1}{p}a$.
\end{proposition}

We shall also need Theorem D from \cite{mpj}.
\begin{theorem}\label{mpj}
There exist absolute constants $c_4,c_5,c_6$ for which the following holds: Let $\mu$ be a probability measure on $\mathbb{R}^n$, $X_1,\ldots,X_m$ be independent random vectors distributed according to $\mu$, set $S$ to be a subset of the unit sphere of $L_2(\mu)$ and assume that $diam(S,\|\cdot\|_{\psi_2})\leq\beta$. Then for any $\theta>0$ and $m\geq1$ satisfying the inequality
\[
  c_4\beta\gamma_2(S,\|\cdot\|_{\psi_2})\leq\theta\sqrt{m}
\]
we have
\begin{equation}\label{mpj1}
  Pr\left(\sup_{f\in S}\left|\frac{1}{m}\sum_{i=1}^{m} f(X_i)^2-1\right|\leq\theta\right)\geq
  1-\exp\left(-\frac{c_5}{\beta^4}\theta^2m\right).
\end{equation}
Moreover, if $S$ is symmetric then
\begin{equation}\label{mpj2}
  \mathbb{E}\sup_{f\in S}\left|\frac{1}{m}\sum_{i=1}^{m}f(X_i)^2-1\right|\leq c_6\max\left\{\beta\frac{\gamma_2(S,\|\cdot\|_{\psi_2})}{\sqrt{m}},
  \frac{\gamma_2(S,\|\cdot\|_{\psi_2})^2}{m}\right\}.
\end{equation}
\end{theorem}
\begin{remark}\label{remark}
Note that if $\beta\leq\frac{\gamma_2(S,\|\cdot\|_{\psi_2})}{\sqrt{m}}$ then \eqref{mpj2} is
\[
  \mathbb{E}\sup_{f\in S}\left|\frac{1}{m}\sum_{i=1}^{m}f(X_i)^2-1\right|\leq
  c_6\frac{\gamma_2(S,\|\cdot\|_{\psi_2})^2}{m}
\]
and if $\beta\geq\frac{\gamma_2(S,\|\cdot\|_{\psi_2})}{\sqrt{m}}$ then \eqref{mpj2} is
\[
  \mathbb{E}\sup_{f\in S}\left|\frac{1}{m}\sum_{i=1}^{m}f(X_i)^2-1\right|\leq c_6\beta\frac{\gamma_2(S,\|\cdot\|_{\psi_2})}{\sqrt{m}}.
\]
\end{remark}

\section{Large distortion embeddings using random variables}

\begin{definition}\label{def100}
Let $X$ be a random vector in $\mathbb{R}^n$ distributed according to a probability measure $\mu$ on $\mathbb{R}^n$.

We say that $\mu$ or $X$ are isotropic if $\mathbb{E}\langle X,a\rangle^2=\|a\|_2^2$ for all $a\in\mathbb{R}^n$.

We say that $\mu$ or $X$ satisfy a $\psi_2$ condition with constant $\beta$ if $\|\langle X,a\rangle\|_{\psi_2}\leq\beta\|a\|_2$ for all $a\in\mathbb{R}^n$.

We say that $\mu$ or $X$ have a concentration property with constant $\alpha$ if $C_\varepsilon\big(\langle X,a\rangle\big)\leq\alpha\varepsilon$ for all $\varepsilon\geq0$ and $a\in S^{n-1}$.
\end{definition}

Let $m,n\in\mathbb{N}$ and $\alpha,\beta>0$. Throughout the rest of this section $\Gamma_1,\ldots,\Gamma_m$ will be centered $n$-dimensional i.i.d. isotropic random vectors which satisfy a $\psi_2$ condition with constant $\beta$ and have a concentration property with constant $\alpha$, and which are distributed according to $\mu$. We put $\Gamma$ to be the matrix whose rows are $\Gamma_1\ldots,\Gamma_m$.

Note that $\beta\geq1$ since for any $a\in\mathbb{R}^n$
\begin{eqnarray}\label{t1}
  &&\|a\|_2=\sqrt{\mathbb{E}\langle\Gamma_1,a\rangle^2}=
  \inf\{C>0:1+\frac{\mathbb{E}\langle\Gamma_1,a\rangle^2}{C^2}\leq2\}\\
  &\leq&\inf\{C>0:\mathbb{E}\exp\frac{\langle\Gamma_1,a\rangle^2}{C^2}\leq2\}=
  \|\langle\Gamma_1,a\rangle\|_{\psi_2}\leq\beta\|a\|_2.\nonumber
\end{eqnarray}

In addition, $\alpha$ is bounded from below: Choosing $a=e_1$ we obtain
\begin{eqnarray*}
  && C_\varepsilon(\langle  \Gamma_1,a\rangle)=C_\varepsilon(\Gamma_{1,1})\leq\alpha\varepsilon\\
  &&\mathbb{E}\langle \Gamma_1,a\rangle^2=\mathbb{E}|\Gamma_{1,1}|^2=\|a\|_2^2=1.
\end{eqnarray*}
Using this and the Chebyshev's inequality we obtain
\begin{eqnarray*}
  &&Pr(|\Gamma_{1,1}|\leq2)\leq2\alpha\\
  &&Pr(|\Gamma_{1,1}|\leq2)=1-Pr(|\Gamma_{1,1}|>2)\geq1-\frac{\mathbb{E}|\Gamma_{1,1}|^2}{4}=\frac{3}{4}
\end{eqnarray*}
which gives $\alpha\geq\frac{3}{8}$.

\medskip

We now present two standard examples of random vectors which satisfy the three condition of definition \ref{def100}

\noindent\textbf{Example:} The random vector with i.i.d. standard Gaussian entries is clearly centered and isotropic and satisfies a $\psi_2$ condition with constant $\sqrt\frac{8}{3}$ and has a concentration property with constant $\sqrt\frac{2}{\pi}$.

\medskip

\noindent\textbf{Example:} The random vector $X=n^\frac{1}{2}U$ where $U$ is uniformly distributed on $S^{n-1}$ is clearly centered and isotropic and has the same distribution as $n^\frac{1}{2}\frac{g}{\|g\|_2}$ where $g=(g_1,\ldots,g_n)$ and the $g_i$ are i.i.d. standard Gaussian variables. Since the random vector is rotation invariant we may choose $a=e_1$. Then, by Lemma 2.2 in \cite{dg}, for $t>2$,
\begin{eqnarray*}
  &&Pr\left(|\langle n^\frac{1}{2}\frac{g}{\|g\|_2},e_1\rangle|\geq t\right)=
  Pr\left(g_1^2\geq\frac{t^2}{n}\sum_{i=1}^{n}g_i^2\right)\\
  &\leq&\exp\left(\frac{1-t^2+\ln{t^2 }}{2}\right)\leq
  \exp\left(-\frac{t^2-4}{4}\right).
\end{eqnarray*}
Since $\exp\left(-\frac{t^2-4}{4}\right)\geq1$ for $0\leq t\leq2$ then
\[
  Pr\left(|\langle n^\frac{1}{2}\frac{g}{\|g\|_2},e_1\rangle|\geq t\right)\leq
  e\exp\left(-\frac{t^2}{4}\right)
\]
for all $t>0$ and by Proposition \ref{psi2} part 2
\[
  \|\langle n^\frac{1}{2}\frac{g}{\|g\|_2},e_1\rangle\|_{\psi_2}\leq (e+1)^\frac{1}{2}2\leq4
\]
and therefore $X$ satisfies a $\psi_2$ condition with constant 4.

In order to estimate $C_\varepsilon(\langle X,a\rangle)$ it is again sufficient to consider $a=e_1$ and, moreover, to consider $L=0$ only in the definition of the $\varepsilon$-concentration of $\langle X,a\rangle$. Then for $0<t<1$, we use Lemma 2.2 in \cite{dg} and obtain
\begin{eqnarray}\label{ye}
  &&Pr\left(|\langle n^\frac{1}{2}\frac{g}{\|g\|_2},e_1\rangle|\leq t\right)=
  Pr\left(g_1^2\leq\frac{t^2}{n}\sum_{i=1}^{n}g_i^2\right)\\
  &\leq&\exp\left(\frac{1-t^2+\ln{t^2}}{2}\right)\leq2t\nonumber
\end{eqnarray}
and since for $t\geq\frac{1}{2}$ we have $2t\geq1$ then \eqref{ye} holds for all $t>0$. Hence $X$ has the concentration property with constant $2$.
\medskip

The following proposition will help us provide many classes of centered isotropic $n$-dimensional random vectors which satisfy a $\psi_2$ condition and have a concentration property.
\begin{proposition}\label{alphabeta}
Let $X=(X_1,\ldots, X_n)$ be a centered random vector in $\mathbb{R}^n$ distributed according to a probability measure $\mu$ on $\mathbb{R}^n$.\\
1. If the entries of $X$ are uncorrelated then $X$ is isotropic if and only if $\mathbb{E}X_i^2=1$ for all $1\leq i\leq n$.\\
2. If $X$ satisfies a $\psi_2$ condition with constant $\beta$ then $\|X_i\|_{\psi_2}\leq\beta$ for $1\leq i\leq n$. If the entries of $X$ are independent and $\|X_i\|_{\psi_2}\leq\beta$ for $1\leq i\leq n$, then $X$ satisfies a $\psi_2$ condition with constant $4\beta$.\\
3. If $X$ has a concentration property with constant $\alpha$ then $C_\varepsilon(X_i)\leq\alpha\varepsilon$ for all $\varepsilon\geq0$ and $1\leq i\leq n$. If the entries of $X$ are i.i.d. then $C_\varepsilon(X_i)\leq\alpha\varepsilon$ for all $\varepsilon\geq0$ and $1\leq i\leq n$ if and only if $\{X_i\}_{i=1}^{n}$ have a density function which is bounded by $\frac{\alpha}{2}$, and in this case, $X$ has a concentration property with constant $c_1\alpha$.\\
\end{proposition}
\proof 1. Obvious.\\
2. The first implication is shown by choosing $a=e_i$ and the second implication follows from Proposition \ref{psi2} part 5.\\
3. The first implication is shown by choosing $a=e_i$.

Assume $C_\varepsilon(X_i)\leq\alpha\varepsilon$ for all $\varepsilon\geq0$ and $1\leq i\leq n$. Let $F_{X_1}$ be the distribution function of $X_1$. Then for $x>y$ put $t=\frac{x+y}{2}$ and $\varepsilon=\frac{x-y}{2}$. Then
\begin{eqnarray*}
  &&F_{X_1}(x)-F_{X_1}(y)=F_{X_1}(t+\varepsilon)-F_{X_1}(t-\varepsilon)=Pr(|X_1-t|\leq
  \varepsilon)\\
  &\leq&C_\varepsilon(X_1)\leq\alpha\varepsilon=\frac{\alpha}{2}(x-y).
\end{eqnarray*}
Hence $F_{X_1}(x)$ has Lipschitz constant at most $\frac{\alpha}{2}$. The converse is obvious. The last claim follows from Theorem \ref{ob}.
\qed

Proposition \ref{alphabeta} implies that if $\{X_i\}_{i=1}^{n}$ are i.i.d. centered random variables which satisfy $\mathbb{E}X_i^2=1$ and $\|X_i\|_{\psi_2}\leq\beta$ for $1\leq i\leq n$, and have a density function which is bounded by $\alpha$, then $X=(X_1,\ldots,X_n)$ is a centered isotropic $n$-dimensional random vector which satisfies a $\psi_2$ condition with constant $4\beta$ and has a concentration property with constant $2c_1\alpha$. This provides a variety of different classes of random matrices whose rows are independent centered isotropic $n$-dimensional random vectors which satisfy a $\psi_2$ condition and have a concentration property.

We present bounds for $\mathbb{E}(S)$. Recall that the Grassman manifold $\mathcal{G}_{m,k}$ is
the collection of all $k$-dimensional subspaces of $\mathbb{R}^m$
equipped with the metric

\[
  \rho(V,W)=\max_{v\in V\cap S^{m-1}}d(v,W\cap S^{m-1})
\]
where $V,W$ are $k$-dimensional subspaces of $\mathbb{R}^m$ and the metric $d$ is Euclidean.
\begin{proposition}\label{prop1}
1. Let $1\leq k\leq n$ and let $\{W_l\}_{l=1}^{p}$ be linear
subspaces of dimension at most $k$ in $\ell_2^n$. Fix $0\leq r<1$. Put
\[
  S_r=\{x\in S^{n-1}:\ There\ exists\ 1\leq l\leq p\ such\ that\
  d(x,W_l\cap S^{n-1})\leq r\}.
\]

Then
\[
  \mathbb{E}(S_r)=\mathbb{E}\max_{x\in S_r}\sum_{i=1}^{n}x_ig_i\leq
  3\big(\sqrt{\ln{p}}+\sqrt{k}+r(\sqrt{\ln{p}}+\sqrt{n-k})\big)
\]
where the $g_i$ are i.i.d. $N(0,1)$ random variables.\\
2. Let $1\leq k\leq n$ and $0<\delta<\sqrt2$ let $\{W_l\}_{l=1}^{p}$ be linear
subspaces of dimension at least $k$ in $\ell_2^n$ which satisfy $\rho(W_i,W_j)\geq\delta$ for all $1\leq i<j\leq p$. Put $S=\cup_{l=1}^{p}W_l\cap S^{n-1}$. Then
\[
  \mathbb{E}(S)=\mathbb{E}\max_{x\in S}\sum_{i=1}^{n}x_ig_i=\Omega\left(\sqrt{k}+\sqrt\frac{\ln{p}}{\ln\frac{1}{\delta}}\right).
\]
where the $g_i$ are i.i.d. $N(0,1)$ random variables.
\end{proposition}
\proof 1. Proven in Proposition 2.5 in \cite{dg2}.\\
2 will be proven later in Proposition \ref{dg6}. \qed

We now turn to state and prove the main Theorem. Recall that $\Gamma$ is an $m\times n$ random matrix whose rows are i.i.d., centered, isotropic, satisfy a $\psi_2$ condition with constant $\beta$ and have a concentration property with constant $\alpha$.
\begin{theorem}\label{largedg}
There is a positive constant $c(\alpha,\beta)$ such that the following
holds: Given $D\geq c(\alpha,\beta)$ and $1\leq k\leq n$ and $p$ affine subspaces
$\{W_l\}_{l=1}^{p}$ of dimension at most $k$ in $\ell_2^n$ and any
$m\geq 5\left(k+\frac{\ln{p}}{\ln{D}}\right)$, there is a number $L$ such that the probability that for
every $1\leq l\leq p$ and $x,y\in W_l$
\[
  \frac{L}{D}\|x-y\|_2\leq\|\Gamma(x)-\Gamma(y)\|_2\leq L\|x-y\|_2
\]
is at least $1-2D^{-\frac{m}{5}}$.
\end{theorem}
\proof We shall need the following estimate. Since $m\geq 5\left(k+\frac{\ln{p}}{\ln{D}}\right)$ then $\ln{p}\leq\frac{1}{5}m\ln{D}$ and $k\leq\frac{m}{5}$ and therefore
\begin{equation}\label{temp66}
  \sqrt{\ln{p}}+\sqrt{k}\leq\sqrt{m\ln{D}}.
\end{equation}

We may assume $\{W_l\}_{l=1}^{p}$ are subspaces of $\ell_2^n$.

Put $S=\cup_{l=1}^{p}\big(W_l\cap S^{n-1}\big)$ and $E^2=\mathbb{E}\max_{x\in S} \|\Gamma(x)\|_2^2$.  Given $A<B$ to be chosen later, we will provide a lower bound for the probability that
\[
  A\cdot E^2\leq\min_{x\in S}\|\Gamma(x)\|_2^2\leq\max_{x\in S}\|\Gamma(x)\|_2^2\leq B\cdot E^2.
\]

Since in $S^{k-1}$, for every $0<\varepsilon<1$ there is an $\varepsilon$-net $N$ with $|N|\leq(\frac{3}{\varepsilon})^{k}$, it follows that there is a
$\sqrt\frac{A}{B}$\ -net $\{x_h\}_{h\in\Delta}$ of $S$ with $|\Delta|\leq p\cdot \left(\frac{9B}{A}\right)^\frac{k}{2}$.

Any $s\in\mathbb{R}^n$ may be considered as an element of $L_2(\mu)$ in the following manner: $s(x)=\langle x,s\rangle$. Since the rows of $\Gamma$ are isotropic then
\[
  \|s\|_{L_2}^2=\int_{\mathbb{R}^n}|\langle x,s\rangle|^2d\mu(x)=
  \mathbb{E}\langle \Gamma_1,s\rangle^2=\|s\|_2^2
\]
which implies that $S$, considered as a subset of $L_2(\mu)$, is a subset of the unit sphere of $L_2(\mu)$. Moreover, by \eqref{t1}
\begin{equation}\label{55}
diam(S,\|\cdot\|_{\psi_2})\leq2\beta
\end{equation}
and the Lipschitz constant of the map $I:(S,\|\cdot\|_2)\longrightarrow(S,\|\cdot\|_{\psi_2})$ is at most $\beta$ and the Lipschitz constant of $I^{-1}$ is at most $1$. Hence by Lemma \ref{lemma02}
\begin{equation}\label{56}
  \gamma_2(S,\|\cdot\|_2)\leq\gamma_2(S,\|\cdot\|_{\psi_2})\leq
  \beta\gamma_2(S,\|\cdot\|_2).
\end{equation}

Note that for any $x\in S^{n-1}$ we have
\begin{eqnarray*}
  &&\mathbb{E}\|\Gamma(x)\|_2^2=
  \mathbb{E}\sum_{i=1}^{m}\langle\Gamma_i,x\rangle^2=
  \sum_{i=1}^{m}\|x\|_2^2=m.
\end{eqnarray*}
Hence
\begin{equation}\label{temp99}
  m\leq E^2.
\end{equation}

Put $\Lambda=\big\{\max_{x\in S}\|\Gamma(x)\|_2^2\leq B\cdot E^2\big\}$. Using conditional probability we obtain
\begin{eqnarray}\label{eq08}
  &&Pr\left(A\cdot E^2\leq\min_{x\in S}\|\Gamma(x)\|_2^2\leq\max_{x\in S}\|\Gamma(x)\|_2^2\leq B\cdot E^2\right)\nonumber\\
  &=&Pr\left(A\cdot E^2\leq\min_{x\in S}\|\Gamma(x)\|_2^2\Big|
  \Lambda\right)\cdot Pr\left(\Lambda\right)\nonumber\\
  &\geq&Pr\left(\bigcap_{h\in\Delta}\Big\{4\cdot A\cdot E^2\leq\|\Gamma(x_h)\|_2^2
  \Big\}\Big|\Lambda\right)Pr(\Lambda)\\
  &=&Pr\left(\bigcap_{h\in\Delta}\Big\{4\cdot A\cdot E^2\leq\|\Gamma(x_h)\|_2^2
  \Big\}\bigcap\Lambda\right)\nonumber\\
  &\geq&1-Pr\left(\bigcup_{h\in \Delta}\Big\{\|\Gamma(x_h)\|_2^2< 4\cdot A\cdot E^2
  \Big\}\right)-Pr(\Lambda^c)\nonumber.
\end{eqnarray}
The proof will be concluded when we show that
\begin{eqnarray*}
  &&1.\ Pr(\Lambda^c)\leq D^{-m}\\
  &&2.\ Pr\left(\bigcup_{h\in \Delta}\Big\{\|\Gamma(x_h)\|_2^2< 4\cdot A\cdot E^2
  \Big\}\right)\leq D^{-\frac{m}{5}}.
\end{eqnarray*}

\noindent\textbf{1. The proof of $Pr(\Lambda^c)\leq D^{-m}$}\\
Using \eqref{temp99} and putting $B=4\beta^2c_5^{-\frac{1}{2}}(\ln{D})^\frac{1}{2}+1>1$ we obtain
\begin{eqnarray}\label{gamma2}
  &&Pr(\Lambda^c)=Pr\left(\max_{s\in S}\|\Gamma(s)\|_2^2>
  B\cdot E^2\right)\nonumber\\
  &\leq&Pr\left(\max_{s\in S}\|\Gamma(s)\|_2^2+m-m> Bm\right)\nonumber\\
  &=&Pr\left(\max_{s\in S}\sum_{j=1}^{n}s(\Gamma_i)^2-m> Bm-m\right)\\
  &\leq&Pr\left(\max_{s\in S}\left|\sum_{j=1}^{n}s(\Gamma_i)^2-m\right|> Bm-m\right)\nonumber\\
  &=&Pr\left(\max_{s\in S}\left|\frac{1}{m} \sum_{j=1}^{n}s(\Gamma_i)^2-1\right|>B-1\right)\nonumber.
\end{eqnarray}
In order to use $\eqref{mpj1}$ to continue \eqref{gamma2} we shall need to verify that $\theta=B-1=4\beta^2c_5^{-\frac{1}{2}}(\ln{D})^\frac{1}{2}$ satisfies the condition $c_4\beta\gamma_2(S,\|\cdot\|_{\psi_2})\leq\theta\sqrt{m}$: By \eqref{56} it is sufficient to show that $c_4\beta^2\gamma_2(S,\|\cdot\|_2)\leq 4\beta^2c_5^{-\frac{1}{2}}(m\ln{D})^\frac{1}{2}$ and using Theorem \ref{ft} with $G_s=\sum_{i=1}^{n}s_ig_i$, and Proposition \ref{prop1}, noting $r=0$, it is sufficient to prove that $3c_4c_2^{-1}(\sqrt{k}+\sqrt{\ln{p}})\leq 4c_5^{-\frac{1}{2}}(m\ln{D})^\frac{1}{2}$ and using \eqref{temp66} in the same manner we get that if $3c_4c_2^{-1}\leq 4c_5^{-\frac{1}{2}}$ then we are done. If $c_5$ does not satisfy this inequality, we can make it smaller without affecting its role in Theorem \ref{mpj}. Now, we use \eqref{mpj1} and obtain
\begin{equation}\label{gamma1}
  Pr(\Lambda^c)\leq\exp\left(-\frac{c_5(B-1)^2m}{(2\beta)^4}\right)=D^{-m}.
\end{equation}

\noindent\textbf{2. The proof of $Pr\big(\cup_{h\in \Delta}\big\{\|\Gamma(x_h)\|_2^2< 4\cdot A\cdot E^2\big\}\big)\leq D^{-\frac{m}{5}}$}\\
\textbf{Case 1:} Assume $\beta\leq\frac{\gamma_2(S,\|\cdot\|_{\psi_2})}{\sqrt{m}}$. Recalling Remark \ref{remark} and using \eqref{mpj2} with the rows of $\Gamma$ as the independent random vectors and $S$ considered as a subset of the unit sphere of $L_2(\mu)$ gives
\begin{eqnarray*}
  &&E^2=\mathbb{E}\max_{s\in S}\|\Gamma(s)\|_2^2=
  \mathbb{E}\max_{s\in S}\sum_{i=1}^{m}\langle\Gamma_i, s\rangle^2=
  \mathbb{E}\max_{s\in S}\sum_{i=1}^{m}s(\Gamma_i)^2\\
  &=&\mathbb{E}\max_{s\in S}\left(\sum_{i=1}^{m}s(\Gamma_i)^2-m+m\right)\leq m+m\mathbb{E}\max_{s\in S}\left|\frac{1}{m}\sum_{i=1}^{m}s(\Gamma_i)^2-1\right|\\ &\leq&m+c_6\gamma_2(S,\|\cdot\|_{\psi_2})^2.
\end{eqnarray*}
Using \eqref{56} we obtain
\begin{eqnarray*}
  &\leq&m+c_6\beta^2\gamma_2(S,\|\cdot\|_2)^2
\end{eqnarray*}
and Theorem \ref{ft} with $G_s=\sum_{i=1}^{n}s_ig_i$, combined with the given symmetric $S$, gives
\begin{eqnarray*}
  &\leq&m+c_6c_2^{-2}\beta^2\big(\mathbb{E}\max_{s\in S}|
  \sum_{i=1}^{n}s_ig_i|\big)^2=m+c_6c_2^{-2}\beta^2\mathbb{E}(S)^2
\end{eqnarray*}
and finally, using Proposition \ref{prop1}, noting that $r=0$, and \eqref{temp66}, we obtain
\begin{eqnarray*}
  &\leq&m+9c_6c_2^{-2}\beta^2(\sqrt{\ln{p}}+\sqrt{k})^2\leq c_7\beta^2 m\ln{D}.
\end{eqnarray*}
Combining this with \eqref{temp99} we obtain
\begin{equation}\label{estimate1}
  m\leq E^2\leq c_7\beta^2 m\ln{D}.
\end{equation}

We now estimate $Pr\left(\cup_{h\in \Delta}\big\{\|\Gamma(x_h)\|_2^2<4\cdot A\cdot E
\big\}\right)$.  To do this we first estimate $Pr(\|\Gamma(x)\|_2^2<\lambda)$. Since $\Gamma_1,\ldots,\Gamma_m$ are independent and have the concentration property with constant $\alpha$ then $C_\varepsilon(\langle\Gamma_i,x\rangle)\leq\alpha\varepsilon$ for $1\leq i\leq m$ and so by Lemma \ref{lemma01}
\begin{equation}\label{eq09}
  Pr(\|\Gamma x\|_2^2<\lambda m)\leq
  (6\alpha)^m(\lambda)^\frac{m}{2}.
\end{equation}
Using \eqref{estimate1}, \eqref{eq09} and $k\leq\frac{m}{5}$ and $\ln{p}\leq\frac{m\ln{D}}{5}$ and putting $A=\frac{4\beta^2c_5^{-\frac{1}{2}}(\ln{D})^\frac{1}{2}+1}{D^2}=\frac{B}{D^2}$ we obtain
\begin{eqnarray}\label{estimate2}
  &&Pr\left(\bigcup_{h\in \Delta}\big\{\|\Gamma(x_h)\|_2^2<4\cdot A\cdot E^2 \big\}\right)\nonumber\\
  &\leq&\sum_{h\in \Delta}
  Pr\left(\|\Gamma(x_h)\|_2^2<4\cdot A\cdot c_7\cdot\beta^2\cdot m\ln{D}\right)\nonumber\\
  &\leq&p\cdot\left(9\frac{B}{A}\right)^\frac{k}{2}\cdot
  (6\cdot \alpha)^m(4\cdot A\cdot c_7\cdot\beta^2\cdot\ln{D})^\frac{m}{2}\\
  &=&\exp\left(\ln{p}+\frac{k}{2}\ln{(9D^2)}+m\ln{(6\cdot \alpha)}+
  \frac{m}{2}\ln\left(4\cdot A\cdot c_7\cdot\beta^2\cdot\ln{D}\right) \right)\nonumber\\
  &\leq&\exp\left(\ln{p}+2k\ln{D}-m\ln{D}+
  \frac{m}{2}\ln{\left(c_8\cdot\alpha^2\cdot\beta^4\cdot
  (\ln{D})^{\frac{3}{2}}\right)}\right)\nonumber\\
  &=&\exp\left(\ln{p}-
  m\ln{D}\left(1-\frac{2k}{m}-
  \frac{\ln\left(c_8\cdot\alpha^2\cdot\beta^4(\ln{D})^{\frac{3}{2}}\right)}{2\ln{D}}
  \right)\right)\nonumber\\
  &\leq&\exp\left(\ln{p}-\frac{2}{5}m\ln{D}\right)\leq
  \exp\left(-\frac{m}{5}\ln{D}\right)=D^{-\frac{m}{5}}\nonumber
\end{eqnarray}
where the inequality before last follows from a suitable choice of large $c_1(\alpha,\beta)$ for which
\[
  \frac{\ln\left(c_8\cdot\alpha^2\cdot\beta^4(\ln{D})^{\frac{3}{2}}\right)}{2\ln{D}}\leq \frac{1}{5}
\]
whenever $D\geq c_1(\alpha,\beta)$.

\textbf{Case 2:} $\beta\geq\frac{\gamma_2(S,\|\cdot\|_{\psi_2})}{\sqrt{m}}$. Remark \ref{remark} implies that \eqref{estimate1} becomes
\[
  m\leq E^2\leq c_7\beta^2 m\sqrt{\ln{D}}.
\]
Then \eqref{estimate2} becomes
\begin{eqnarray*}
  &&Pr\left(\bigcup_{h\in \Delta}\big\{\|\Gamma(x_h)\|_2^2<4\cdot A\cdot E^2 \big\}\right)\\
  &\leq&\exp\left(\ln{p}-
  m\ln{D}\left(1-\frac{2k}{m}-
  \frac{\ln\left(c_8\cdot\alpha^2\cdot\beta^4\ln{D}\right)}{2\ln{D}}
  \right)\right)\\
  &\leq&\exp\left(\ln{p}-\frac{2}{5}m\ln{D}\right)\leq
  \exp\left(-\frac{m}{5}\ln{D}\right)=D^{-\frac{m}{5}}
\end{eqnarray*}
where the inequality before last follows from a suitable choice of large $c_2(\alpha,\beta)$ for which
\[
  \frac{\ln\left(c_8\cdot\alpha^2\cdot\beta^4\ln{D}\right)}{2\ln{D}}\leq\frac{1}{5}
\]
whenever $D\geq c_2(\alpha,\beta)$.

Setting $c(\alpha,\beta)=c_1(\alpha,\beta)$ when $\beta\leq\frac{\gamma_2(S,\|\cdot\|_{\psi_2})}{\sqrt{m}}$ and $c(\alpha,\beta)=c_2(\alpha,\beta)$ when $\beta\geq\frac{\gamma_2(S,\|\cdot\|_{\psi_2})}{\sqrt{m}}$ concludes the proof. Note that $L=E\sqrt{B}$.
\qed

\begin{remark}
The constant $c(\alpha,\beta)$ can be estimated by
\begin{eqnarray*}
    c(\alpha,\beta)\leq c_9\big(\alpha^2\beta^4\big)^\frac{5}{2}\big(\ln(1+\alpha^2\beta^4)\big)^{\frac{15}{4}}
\end{eqnarray*}
\end{remark}

\begin{remark}
The estimate on the dimension in Theorem \ref{largedg} is tight in terms of $k,p,D$ as was proved in Theorem 1.2 part (iii) in \cite{dg2}.
\end{remark}

\begin{remark}
Note that if we wish to take $x\in W_l$ and $y\in W_{l'}$ and still distort the distance between them by a factor of $D$ at most, then the estimate for $m$ barely changes because we need only take $V_{l,l'}=span\{W_{l}'\cup W_{l'}'\}$ where $W_l'$ is the unique parallel subspace to $W_l$. Apply Theorem \ref{largedg} to $\{V_{l,l'}\}_{l,l'=1}^{p}$ and since $dim(V_{l,l'})\leq 2k$ and the number of subspaces is not more than $p^2$, the estimate on $m$ is larger than the estimate for $\{W_l\}_{l=1}^{p}$ by a factor of $2$.
\end{remark}

The proof of Theorem \ref{largedg} can be modified slightly to prove the following:
\begin{theorem}\label{largedg2}
There is a positive constant $c(\alpha,\beta)$ such that the following
holds: Given $D\geq 2$ and $1\leq k\leq n$ and $p$ affine subspaces
$\{W_l\}_{l=1}^{p}$ of dimension at most $k$ in $\ell_2^n$ and any
$m\geq c(\alpha,\beta)\left(k+\frac{\ln{p}}{\ln{D}}\right)$, there is a number $L$ such that the probability that for
every $1\leq l\leq p$ and $x,y\in W_l$
\[
  \frac{L}{D}\|x-y\|_2\leq\|\Gamma(x)-\Gamma(y)\|_2\leq L\|x-y\|_2
\]
is at least $1-2D^{-\frac{m}{5}}$.
\end{theorem}

\begin{remark}
The constant $c(\alpha,\beta)$ in theorem \ref{largedg2} is $c\alpha^2\beta^4$ where $c$ is a universal constant.
\end{remark}

The proof of part 2 of Proposition \ref{prop1} follows from the following proposition.
\begin{proposition}\label{dg6}
Let $1\leq k\leq m\leq n$ and $0<\delta\leq\sqrt2$ and $D\geq2$ and let $\{W_l\}_{l=1}^p$ be subspaces of dimension at least $k$ in $\ell_2^n$ and put $S=\cup_{l=1}^{p}(W_l\cap S^{n-1})$.\\
If $\rho(W_l,W_{l'})\geq\delta$ for all $1\leq l<l'\leq p$, then any linear map $F:\ell_2^n\longrightarrow\ell_2^m$ which distorts the Euclidean norm of the elements of $S$ by factor $D$ at most satisfies $m=\Omega\big(k+\frac{\ln{p}}{\ln\frac{D}{\delta}}\big)$. In addition $\mathbb{E}(S)=\Omega\big(\sqrt{k}+ \sqrt\frac{\ln{p}}{\ln\frac{1}{\delta}}\big)$.
\end{proposition}
\proof

For $A,B\subset\mathbb{R}^n$ denote by $d_{aff}(A,B)$ the distance between the
affine subspace spanned by $A$ to the affine subspace spanned by $B$. If $A=\{x\}$ then we write $d_{aff}(x,B)$.

Let $D\geq 2$. We may assume $k\leq \ln{p}$. For every $1\leq i\neq j\leq p$ choose $v_{i,j}\in W_i\cap S^{n-1}$ such that $d(v_{i,j},W_j\cap S^{n-1})=\rho(W_i,W_j)\geq\delta$. Then $d(v_{i,j},W_j)\geq\frac{\delta}{2}$. Let $B_j=\{b_{j,0}=0,b_{j,1},\ldots,b_{j,k}\}$ where $\{b_{j,1},\ldots,b_{j,k}\}$ is an independent set in $W_j$. Then
\begin{equation}\label{mama}
  \frac{\delta}{2}\leq d(v_{i,j},W_j)=d_{aff}(v_{i,j},B_j).
\end{equation}

We first prove the lower bound on $m$. Assume
$F:\ell_2^n\longrightarrow\ell_2^m$ is any map which satisfies $F(0)=0$ and
such that for every $1\leq i\neq j\leq p$
\begin{eqnarray}\label{eq900}
  &&\frac{1}{D}\leq
  \frac{d_{aff}(F(v_{i,j}),F(B_j))}{d_{aff}(v_{i,j},B_j)}\leq 1\\
  &&\frac{1}{D}\leq\frac{d_{aff}(F(v_{i,j}),0)}{d_{aff}(v_{i,j},0)}=
  \frac{d(F(v_{i,j}),0)}{d(v_{i,j},0)}=\|F(v_{i,j})\|_2\leq 1\nonumber
\end{eqnarray}
and for every $1\leq i\leq p$ and $1\leq l\leq k$
\begin{equation}\label{eq901}
  \frac{1}{D}\leq
  \frac{d_{aff}(F(b_{i,l}),F(\{b_{i,0},\ldots ,b_{i,l-1}\}))}
  {d_{aff}(b_{i,l},\{b_{i,0},\ldots ,b_{i,l-1}\})}\leq 1.
\end{equation}

Put $V_i=span(F(B_i))$. Then by \eqref{eq901} $\{V_i\}_{i=1}^{p}$ are
$k$-dimensional subspaces of $\ell_2^m$. Consider $V_i$ as elements
of the Grassman manifold $\mathcal{G}_{m,k}$. Then by \eqref{mama} and \eqref{eq900} $\{V_i\}_{i=1}^{p}$
satisfy $\rho(V_i,V_j)\geq\frac{\delta}{2D}$ for every $1\leq i\neq j\leq
p$. Hence the Grassman manifold $\mathcal{G}_{m,k}$ contains $p$
disjoint balls of radius $\frac{\delta}{4D}$ and therefore
$p\cdot\mu_{m,k}(B_\rho\big(V_1,\frac{\delta}{4D})\big)\leq1$ where $\mu_{m,k}$ is the normalized Haar measure on
$\mathcal{G}_{m,k}$. It is well known that there is a universal constant $C>0$ such that for every
$0<r<\frac{1}{2}$
\[
  \mu_{m,k}\big(B_\rho(W,r)\big)\geq (Cr)^{m-k}.
\]
Hence $m\geq c(k+\frac{\ln{p}}{\ln\frac{D}{\delta}})$ for
some universal constant $c>0$.
We now observe that if $F:\ell_2^n\rightarrow\ell_2^m$ is linear and
for all $x\in S$
\[
  \frac{1}{D}\leq\|F(x)\|_2\leq1
\]
then $F$ satisfies \eqref{eq900} and
\eqref{eq901} and therefore $m\geq c(k+\frac{\ln{p}}{\ln\frac{D}{\delta}})$.

Now, Using Theorem 2.3 part i in \cite{dg2} with $\epsilon=\frac{1}{2}$ we obtain a linear map $F:\ell_2^n\longrightarrow\ell_2^m$ with $m=8\mathbb{E}(S)^2$ and $D=2$. Hence $8\mathbb{E}(S)^2\geq c(k+\frac{\ln{p}}{\ln\frac{1}{\delta}})$ and therefore $\mathbb{E}(S)\geq c'(\sqrt{k}+\sqrt\frac{\ln{p}}{\ln\frac{1}{\delta}})$. \qed

\begin{remark}
Let $A=\{e_1,\ldots,e_n\}\subset\ell_2^n$ and $\{W_l\}_{1\le l\le {n \choose k}}$ be the collection of all subspaces made up of all spans of $k$-sparse vectors, i.e. the subspaces are $span\{e_{i_1},\ldots,e_{i_k}\}$ where $i_1<\ldots<i_k$. We obtain ${n \choose k}$ subspaces for which $\rho(W_i,W_j)=\sqrt{2}$ whenever $1\leq i\neq j\leq {n \choose k}$. Hence if $1\le p\le {n \choose k}$,  taking any $p$ of the ${n \choose k}$ subspaces and putting $S=\cup_{l=1}^p (W_l\cap S^{n-1})$ gives $\mathbb{E}(S)\sim \sqrt{k}+\sqrt{\ln{p}}$ and the corresponding lower estimate on the dimension $m$ in Proposition \ref{dg6} is $m\ge c(k+\frac{\ln{p}}{\ln{D}})$.
\end{remark}

\section{Appendix}

\textbf{Proof of Proposition \ref{psi2}:}\\
1. Let $h>0$. Using the Chebychev inequality we obtain
\begin{eqnarray*}
  Pr(|X|>t)=Pr(\exp(h|X|^p)>\exp(ht^p))\leq\mathbb{E}\exp(h|X|^p)\exp(-ht^p).
\end{eqnarray*}
Let $h$ tend to $\|X\|_{\psi_p}^{-p}$ from below and use the definition of the $\psi_p$ norm to obtain the desired result.\\
2. Let $C>a$. Since $\exp\frac{|X|^p}{C^p}$ is positive, then
\begin{eqnarray*}
  &&\mathbb{E}\exp\frac{|X|^p}{C^p}=\int_0^\infty Pr\left(\exp\frac{|X|^p}{C^p}>t\right)dt\\
  &=&1+\int_1^\infty Pr\left(|X|>C(\ln{t})^\frac{1}{p}\right)dt\leq
  1+\int_1^\infty b\cdot\exp\left(-\frac{C^p\ln{t}}{a^p}\right)dt\\
  &=&1+b\int_1^\infty t^{-\frac{C^p}{a^p}}dt=1+\frac{b}{\frac{C^p}{a^p}-1}\leq2
\end{eqnarray*}
where the last inequality holds whenever $C\geq(b+1)^{\frac{1}{p}}a$ which implies $\|X\|_{\psi_p}\leq(b+1)^\frac{1}{p}\cdot a$.\\
3. For $h>0$,
\begin{eqnarray}\label{eq10}
  &&Pr(X>t)=Pr(\exp(hX)>\exp(ht))\nonumber\\
  &\leq&\mathbb{E}\exp(hX)\exp(-ht)\leq\exp(c^2h^2-ht).
\end{eqnarray}
Minimizing the expression on $h>0$ we obtain $2c^2h-t=0$ or $h=\frac{t}{2c^2}$. Substituting into \eqref{eq10} we obtain
\[
  Pr(X>t)\leq\exp\left(-\frac{t^2}{4c^2}\right).
\]
Similarly
\[
  Pr(X<-t)\leq\exp\left(-\frac{t^2}{4c^2}\right)
\]
hence
\[
  Pr(|X|>t)\leq2\exp\left(-\frac{t^2}{4c^2}\right)
\]
and therefore by part 2 we obtain the desired result.\\
4. Let $C>a$. Then
\begin{eqnarray*}
  2\geq\mathbb{E}\exp\frac{X^2}{C^2}=\sum_{n=0}^\infty\frac{\mathbb{E}X^{2n}}{C^{2n}n!}
\end{eqnarray*}
and therefore $\mathbb{E}X^{2n}\leq2C^{2n}n!$. Hence $\mathbb{E}|X|^{2n-1}\leq2C^{2n-1}n!$ for $n\geq2$. Using this and the simple inequalities $\frac{2n!}{(2n)!}\leq\frac{1}{n!}$ and $\frac{2(n+2)!}{(2n+3)!}\leq\frac{1}{n!}$ for $n\geq1$ we obtain that for $t\geq0$
\begin{eqnarray}\label{ma}
  &&\mathbb{E}\exp(tX)=\sum_{n=0}^\infty\frac{t^n\mathbb{E}X^n}{n!}=
  1+\sum_{n=1}^\infty\frac{t^{2n}\mathbb{E}X^{2n}}{(2n)!}+
  \sum_{n=2}^\infty\frac{t^{2n-1}\mathbb{E}X^{2n-1}}{(2n-1)!}\nonumber\\
  &\leq&1+\sum_{n=1}^\infty\frac{2(C^2t^2)^{n}n!}{(2n)!}+
  \sum_{n=2}^\infty\frac{2C^{2n-1}t^{2n-1}n!}{(2n-1)!}\nonumber\\
  &\leq&1+\sum_{n=1}^\infty\frac{(C^2t^2)^{n}}{n!}+
  C^3t^3\sum_{n=0}^\infty\frac{2(C^2t^2)^n (n+2)!}{(2n+3)!}\\
  &\leq&\exp(C^2t^2)+C^3t^3\sum_{n=0}^\infty\frac{(C^2t^2)^n}{n!}=
  (1+C^3t^3)\exp(C^2t^2)\leq\exp(2C^2t^2).\nonumber
\end{eqnarray}
Since $\|-X\|_{\psi_2}=1$ and $\mathbb{E}(-X)=0$ then $\mathbb{E}\exp(t(-X))\leq\exp(2C^2t^2)$ for $t\geq0$ and therefore $\mathbb{E}\exp(tX)\leq\exp(2C^2t^2)$ for all $t$ . Since this holds for all $C>a$ then $\mathbb{E}\exp(tX)\leq\exp(2a^2t^2)$.\\
The proof in the case where the random variables are symmetric is the same noting that $\mathbb{E}X^{2n-1}=0$ for $n\geq1$ and therefore the third summand in \eqref{ma} disappears.\\
5. Using part 4 we obtain that $\mathbb{E}\exp(tX_i)\leq\exp(2a^2t^2)$ for $1\leq i\leq n$. Since the random variables are independent then
\begin{eqnarray*}
  \mathbb{E}\exp(t\sum_{i=1}^{n}a_iX_i)=\prod_{i=1}^{n}\mathbb{E}\exp(a_itX_i)\leq
  \prod_{i=1}^{n}\exp(2a^2a_i^2t^2)=\exp(2a^2t^2).
\end{eqnarray*}
Using part 3 we conclude that $\|\sum_{i=1}^{n}a_iX_i\|_{\psi_2}\leq4a$. Similarly when the random variables are symmetric.\\
6.
\begin{eqnarray*}
  \mathbb{E}\exp(t\sum_{i=1}^{n}a_iX_i)=\prod_{i=1}^{n}\mathbb{E}\exp(a_itX_i)\leq
  \prod_{i=1}^{n}\exp(c^2a_i^2t^2)=\exp(c^2t^2).
\end{eqnarray*}
7. Since $Pr(|X|>t)=0$ for any $t\geq a$ and since for $0<t<a$ and $b>0$ we have $\exp\left(-\frac{t^p}{b^p}+\frac{a^p}{b^p}\right)\geq1$ then
\begin{eqnarray*}
  Pr(|X|>t)\leq\exp\left(\frac{a^p}{b^p}\right)\exp\left(-\frac{t^p}{b^p}\right)
\end{eqnarray*}
for all $t\geq0$ and therefore by part 2 $\|X\|_{\psi_p}\leq(\exp\left(\frac{a^p}{b^p}\right)+1)^\frac{1}{p}b$. Substituting $b=a$ gives the desired result.
\qed

\textit{Acknowledgements}: The authors thank Krzysztof Oleszkiewicz for bringing to our attention the result by Ball and Nazarov.

\smallskip

\noindent
A. Dmitriyuk, {\small Dept. of Math.}, {\small Technion},
{\small Haifa 32000, Israel}, {\small \tt alondm@gmail.com}

\smallskip

\noindent
Y. Gordon, {\small Dept. of Math.}, {\small Technion},
{\small Haifa 32000, Israel}, {\small \tt gordon@techunix.technion.ac.il}


\begin{thebibliography}{}


%\bibitem{noga} Alon, N., \emph{Problems and results in extremal
%combinatorics, Part 1}, Unpublished manuscript.

\bibitem{achlioptas} D. Achlioptas,
\emph{Database-friendly random projections: Johnson-Lindenstrauss
with binary coins}, Special issue on PODS 2001 (Santa Barbara,
CA), vol. 66, 2003, 671--687.

\bibitem{av} R.\ I. Arriaga and S. Vempala,
\emph{An algorithmic theory of learning: robust concepts and random
projection}, 40th Annual Symposium on Foundations of Computer
Science (New York, 1999), 1999, 616--623.

\bibitem{bo} K. Ball, and F. Nazarov, \emph{Little level theorem and Zero-Khinchin inequality for sums of independent random variables}, 1996, http://www.mth.msu.edu/~fedja/prepr.html

%\bibitem{bk} F. Barthe and A. Koldobski,
%\emph{Extremal slabs in the cube and the Laplace transform}, Adv.
%Math., vol. 174, No. 1, 2003, 89--114.

%\bibitem{bs} J. Bourgain and S. J. Szarek,
%\emph{The Banach-Mazur distance to the cube and the Dvoretzky-Rogers
%factorization}, Israel J. Math., vol. 62, 1988, 169--180.

\bibitem{bourgain} Bourgain, J.,
\emph{On {L}ipschitz embedding of finite metric spaces in {H}ilbert
              space}, Israel J. Math., Vol. 52, 1985, 46--52.

\bibitem{ct} E. J. Candes and T. Tao,
\emph{Near-optimal signal recovery from random projections:
universal encoding strategies}, IEEE Trans. Inform. Theory, 52,
2006, 5406--5425.

\bibitem{dg} S. Dasgupta and A. Gupta, \emph{An
elementary proof of a theorem of Johnson and Lindenstrauss},
Random Structures Algorithms, vol. 22, 2003, 60--65.

\bibitem{dg2} A. Dmitriyuk and Y. Gordon, \emph{Generalizing the Johnson-Lindenstrauss lemma to k-dimensional affine subspaces},
Studia Math., 195(2009), 2009, 227--241.

\bibitem{donoho} D.L Donoho, \emph{Compressed sensing},
IEEE Trans. Inform. Theory, vol. 52, 2006, 1289--1306.

\bibitem{fernique} X. Fernique, \emph{R\'egularit\'e des trajectoires des fonctions aleatoires gaussiennes},
Ecole d'Et\'e de Probabilit\'es de Saint-Flour IV-1974, Lecture Notes in Mathematics 480, Springer-Verlag 1975, 1--96.

%\bibitem{figiel} Figiel, T. and Tomczak-Jaegermann, Nicole, \emph{
%Projections onto {H}ilbertian subspaces of {B}anach spaces}, Israel J. Math., vol. 33, %1979, 155--171.

%\bibitem{fm} P. Frankl and H. Maehara, \emph{The
%Johnson-Lindenstrauss lemma and the sphericity of some graphs}, J.
%Combin. Theory Ser. B, vol. 44, 1988, 355--362.

%\bibitem{gluskin} E. D. Gluskin,
%\emph{An octahedron is poorly approximated by random subspaces},
%Funktsional. Anal. i Prilozhen., vol. 20, 1986, 14--20, 96.

\bibitem{gordon85} Y. Gordon, \emph{Some
inequalities for Gaussian processes and applications}, Israel J.
Math., vol. 50, 1985, 265--289.

\bibitem{gordon87} Y. Gordon, \emph{
Elliptically contoured distributions}, Probab. Theory Related
Fields, vol. 76, 1987, 429--438.

\bibitem{gordon88} Y. Gordon, \emph{
On Milman's inequality and random subspaces which escape through a
mesh in {${\bf R}\sp n$}}, Geometric aspects of functional
analysis (1986/87), Lecture Notes in Math., vol. 1317, 1988,
84--106.

\bibitem{ggmp} Y. Gordon, O. Gu{\'e}don, M. Meyer and A.
Pajor, \emph{Random Euclidean sections of some classical Banach
spaces}, Math. Scand., vol. 91, 2002, 247--268.

\bibitem{glsw} Y. Gordon, A. Litvak, C. Sch\"{u}tt, E. Werner, \emph{On the
minimum of several random variables}, Proc. AMS 134 (2006), No. 12,
3665-3675.


\bibitem{im} P. Indyk and R. Motwani,
\emph{Approximate nearest neighbors: towards removing the curse of
dimensionality}, STOC '98 (Dallas, TX), 1999, 604--613.

\bibitem{jl} W.\ B. Johnson. and J. Lindenstrauss,
\emph{Extensions of Lipschitz mappings into a Hilbert space},
Conference in Modern Analysis and Probability (New Haven, Conn.,
1982), Contemp. Math., vol 26, 1984, 189--206.

\bibitem{magen} A. Magen, \emph{Dimensionality reductions
that preserve volumes and distance to affine spaces, and their
algorithmic applications}, Randomization and approximation
techniques in computer science, Lecture Notes in Comput. Sci., vol.
2483, 239--253.

\bibitem{mpj} S. Mendelson, A. Pajor and N. Tomczak-Jaegermann,
\emph{Reconstruction and subgaussian operators in asymptotic geometric analysis},
Geom. Funct. Anal., vol. 17, 1248--1282.

%\bibitem{mp} Gilles Pisier, \emph{Probabilistic
%methods in the geometry of Banach spaces}, Probability and Analysis
%(Varenna, 1985), Lecture Notes in Math., vol. 1206, 1986, 167--241.

\bibitem{Talagrand} M. Talagrand, \emph{Regularity of Gaussian processes},
Acta Math. 159 (1987), 99--149.
\end{thebibliography}
\end{document}